\documentclass[hidelinks,onefignum,onetabnum]{siamart250211} 

\usepackage[margin=1in,letterpaper]{geometry}
\usepackage{lipsum}
\usepackage{amsfonts}
\usepackage{graphicx}
\usepackage{epstopdf}
\usepackage{algorithmic}
\ifpdf
  \DeclareGraphicsExtensions{.eps,.pdf,.png,.jpg}
\else
  \DeclareGraphicsExtensions{.eps}
\fi

\newsiamremark{remark}{Remark}
\newsiamremark{hypothesis}{Hypothesis}
\crefname{hypothesis}{Hypothesis}{Hypotheses}
\newsiamthm{claim}{Claim}
\newsiamremark{fact}{Fact}
\crefname{fact}{Fact}{Facts}
\newcommand{\ds}{\displaystyle}
\headers{ROM and Anderson type acceleration for DNN}{K. Ito and T. Xue}

\title{Anderson-type acceleration method for Deep Neural Network optimization \thanks{Submitted to the editors \today.
}} 

\author{Kazufumi Ito\thanks{Department of Mathematics, Graduate Program of Operations Research, Center for Research in Scientific Computation, North Carolina State University, Corresponding author.
  (\email{kito@ncsu.edu}, \url{https://kito.wordpress.ncsu.edu/}).}
\and Tiancheng Xue \thanks{Grauate Program of Operations Research, North Carolina State University
  (\email{txue2@ncsu.edu}).}
} 

\usepackage{amsopn}

\newcommand{\R}{\mathbf{R}}
\newcommand{\Z}{\mathbf{Z}}

\ifpdf
\hypersetup{
  pdftitle={Anderson},
  pdfauthor={Kazufumi Ito， Tiancheng Xue}
}
\fi

\begin{document}

\maketitle

\begin{abstract}
In this paper we consider the neural network optimization. We develop Anderson-type acceleration method for the stochastic gradient decent method and it improves the network permanence very much. We demonstrate the applicability of the method for Deep Neural Network (DNN) and Convolution Neural Network (CNN). 

\end{abstract}



\section{Introduction}\label{intro}
In this paper we consider Anderson-type acceleration method \cite{A}  \cite{IX} \cite{WN} for neural network optimization.  
Neural network \cite{MP} is named and motivated by neural networks to mimic the human brain by using interconnected nodes called neurons to process data and learn patterns. They  create connections and adjust the weights of those connections through a process called training, which involves feeding the network large datasets. This training allows the network to refine its weights to minimize errors between its predictions and actual values, similar to how humans learn by recognizing patterns. 
Neural network optimization is the process of tuning a network's parameters (weights and biases) to minimize a loss function, which quantifies the difference between predicted and actual values, ultimately improves the model's performance. Key techniques include using optimization algorithms like Gradient Descent to iteratively tune hyper-parameters. In particular, one may choose appropriate optimizers (e.g., Adam, SGD) \cite{GBC} \cite{KB}.  
We use Neural Network to formulate the forward path and category (digit) classification criterion as the performance, it results in the neural network optimization.  In our paper, the stochastic gradient decent method is applied as the learning method. 

Anderson-type acceleration (AA) method accelerates the 
learning and results in accurate weights in the network design. Anderson acceleration, also called Anderson mixing, is a method for the acceleration of the convergence rate of fixed-point iterations. 
It works by combining the most recent iterates and update steps in a fixed-point iteration to improve the convergence properties of the sequence. AA is particularly effective in scenarios where the fixed-point iteration might be slow or unstable, leading to faster convergence and improved solution quality. 

We develop the Anderson-type  method to accelerate the learning \cite{RHW} of a given neural network. The basic idea is that it uses the optimal assembling, bootstrapping and extrapolation by the optimal combination of the sequence of network designs \cite{A}.
We minimize the weighted total sum of the residuals corresponding  to the sequence of designs over the weights.  We generate the sequence of designs 
for randomly sorted epochs sequentially. We then apply the Anderson-type acceleration. It should decrease the merit function which defines the performance of  the linear combination of the sequence of designs based on the reduced order method (ROM). Like GMRES, we restart the algorithm with the Anderson update to improve the performance. 

In this paper we develop Anderson-type acceleration method for deep neural network (DNN) \cite{HSW}, i.e.,
we build the  neural network transformation $\psi(W,b)$  of input vectors $A\in \R^{n_0 \times N}$ to output data $Y\in \R^{n_L \times N}$ by the discrete time dynamics 
\begin{equation} \label{DNN}
\begin{array}{l}X_{k+1}=\phi_k(W_kX_k+b_k),\;\; X_0=A
\\ \\
Y=H(X_L)
\end{array} \end{equation}
where the labels $X_k\in \R^{n_k\times N}$ and $N$ is the sample size, $\phi_k$ is the activation function such as Relu $\phi(x) = \max(0,x)$ for layer $k$ ($1\le k < L$) and $H$ is the target constraint map \cite{Ho}.
Here $L$ is the number of the layers and
the weight matrices $W_k\in \R^{n_{k-1}\times n_{k}}$ and the bias vectors $\R^{n_k}$ are the internal design variables to be determined.

In the case the category (digit) classification \cite{LC}, we vectorize the digit as follows:
In the context of computer science and machine learning, "one-hot" refers to a method of representing categorical data as binary vectors $y$ where only one element in the vector is ``hot" (set to 1), indicating the presence of a specific category, and all other elements are ``cold" (set to 0). This technique is used to convert categorical variables into a numerical format suitable for use in machine learning algorithms that require numerical input.  

In this case $H$ is the softmax transforms  $X_L$ to the binary distribution and we then formulate the least squares problem:
\begin{equation}\label{SoMax}
\mathcal{L}(y,x_L) = \frac{1}{2N}|H(x_L)-Y|^2, \text{and} \ H(z) = \frac{1}{\sum_{j=1}^{n_L} e^{z_j}} (e^{z_1}, \cdots, e^{z_{n_L}}),    
\end{equation}
where we vectorize the actual result $y$ with ``one-hot'' to $Y$: We represent categorical data as binary vectors $Y$ where only one element in the vector is ``hot'' (set to 1), indicating the presence of a specific category, and all other elements are ``cold'' (set to 0). 

In general  the model-based transport map $T_\theta$, one can formulate the least squares problem for the transport problem:
given the data set $(g_0,g_1)$
\begin{equation} \label{opt}
\min\quad |T_\theta(g_0)-g_1|^2+\Theta(\theta) \mbox{  over $\theta$} 
\end{equation}
where $T_{\theta}$ is the transport map and 
$\Theta$ is the transport cost. The map can represent the tomographic map, e.g. CT \cite{AMOS}, EIT  \cite{IJ}, acoustic, imaging \cite{GA}, electro-magnetic scattering, the Monge-Katrovitch mass transport and the value and mass transport by PDE based models \cite{Itma} with $\theta$ as the corresponding medium parameter. We have DNN modeling of the input-output map $T_\theta$ using pairs of input  and out put data by means of the universal approximation theory \cite{HSW} \cite{C}. Thus, DNN is a particular case of the transportation problem but is a model free method.  The method  and algorithms  developed in the paper can be applied to the general case. We discuss the Anderson-type acceleration method to determine $\theta$ based on DNN model in Section \ref{Andmethod}. 

Another important class of Neural Network is Convolution Neural Network. Convolution neural-network (CNN) \cite{LBBH} \cite{KSH} transfers images from  $A \to B$ by convolution and then other operations (such as max pooling) with a fully connected layer $B\to Y$,
where $B$ has much richer information than  $A$ and $B\to Y$ is performed by fully-connected layer. As we shall see, a fully-connected layer is a single layer DNN. The data compression and filtered decomposition are the key aspect of convolution operation. Other model based network such as DEnet \cite{SP}, Voltera \cite{MZ} net and mass transportation can be formulated similar to DNN as forward-path models, e.g, forecast and prediction are by the models.

Before introducing paper structure, we list some applications of the Neural Network formulation and Anderson-type acceleration. 
\begin{itemize}
    \item Model Informed method (e.g., PINN) \cite{RPK} \cite{AMOS} \cite{RPK} \cite{DeepONet2021} uses the neural network basis elements with the internal design variables to represent solutions and parameters and then the least squares method to equate the model equations.
    \item In \cite{IRZ} we apply to solve the Hamilton-Jacobi equation for the game theory we accurately compute the feedback law via the neural network function.
\end{itemize}

\noindent \textbf{\underline{Contribution of the paper:}} The paper extends the range of applications developed in \cite{IX} with a randomized initialization, and increased problem ill-posedness. 

\noindent \textbf{\underline{Paper Outline}}: We first give highlights of ROM, and Anderson type acceleration developed in \cite{IX}.
We then introduce essential components for a Neural Network, including forward propagation, backward propagation, the standard stochastic gradient decent method. We introduce our numerical testing at the end of the paper. 

\noindent \textbf{\underline{Notation:}} Throughout the paper, we use $x^+$ to denote solution update.

\section{ROM, Anderson type acceleration}
We give some highlighting arguments from \cite{IX}. 
\subsection{Anderson type acceleration method}\label{Andmethod}
Consider the least squares problem
$$
|\psi(x)-y|^2,
$$
where $\psi:X\to Y$ is a semi-smooth map.
\begin{definition}[Anderson type method]
Let $\{x_i\}_{i=1}^m$ be a sequence of approximated solutions to the least square  problem. The modified Anderson method minimizes the total residual
\begin{equation} \label{And}
\min_{\alpha}\ \sum_{i=1}^m |\alpha_i(\psi(x_i)-y)|^2
\end{equation}
subject  to
$$
\sum_{i=1}^m \alpha_i =1.
$$
Then, with optimal solution $\alpha^*$, update $x_{m+1}$ by
$$
x^+=\sum_{i=1}^m \alpha_i^*x_i.
$$    
\end{definition}


\subsection{Reduced order method (ROM)} The Reduced order method (ROM) \cite{IR} gives the theoretical
foundation of the Anderson-type acceleration method \cite{IX}.
\begin{definition}[Reduced order method]
Given $F: X \to X$ in residual form, ROM step gives us the convergence based on the criterion:  $\alpha^*_k$ minimizes
\begin{equation} \label{ROM}
\min\quad |F(\sum_{k\le m}\alpha_k\,x_k)|^2, 
\end{equation}
and update by
$$
x_{m+1}=\sum_{k\le m} \alpha^*_k x_k.
$$    
\end{definition} 

In general let
$$
F(x)=\psi(x)-y
$$
Consider the alternative method based on reduced order
method (ROM) \cite{IR} merit function
\begin{equation} \label{ROM}
\min_{\alpha} \ |\psi(\sum_{k=1}^m\alpha_ix_k)-y|^2.
\end{equation}
It is theoretically advantagous that this guarantees that the merit decreases as the iterate proceeds.  But, the Anderson-type is less expensive to perform.
\subsection{Comparison of Anderson and ROM}
Next, we compare the Anderson and the ROM method.
\begin{theorem}
    Let $F: X\to X$, consider 2 problems: the first one is \eqref{And}, i.e., with $r_k = \psi(x_k) - y$, we have: 
$$
\min\quad|\sum_{k\le m}\alpha_k\,r_k|^2 \mbox{  subject to  } \sum_{k \le m} \alpha_k=1.
$$ 
and \eqref{ROM} with the constraint $\sum_{k\le m}\alpha_k=1$. Then, for some function $H: X\to X$ such that $\ds \lim_{x_k \to \bar{x}} H(x_k) = 0$ for every $\bar{x}$ close to $x_k$, we have:
$$ \sum_k \alpha_k r_k \sim F(\sum_k \alpha_k x_k) + H(x_k),$$
\end{theorem}

\section{Data compression and CNN}
For the fully-connected network the number of parameterizations 
$(W_k,b_k)$
can be very large and expensive. In order to reduce the complicity without loosing the performance, we may constrain weight $W$ by the sparsity and bandwidth by $\ell^0$ sparsity penalization \cite{IK} to constrain the number of nonzero elements in $W$: 
$$
\min\quad |H(X_L)-Y|^2 +\beta |W|_0.
$$





\subsection{Convolution  CNN,  data compression}
In this section we discuss  the convolution network CNN (convolution neural network).
First, we apply the zero-padding with zeros around the border of 
2-d image $X$.
The convolution $X^+\in R^{n,n,c}$ of $X\in R^{n,n}$ with the weight $W$ of size $B=b\times b$ is given by
\begin{equation}\label{conv}
X^+_{i,j,k}=\sum_{(i',j')\in B}W_{i',j',k}X_{i+i',j+j'},
\end{equation}
where  $B$ is a window with size $b\times b$  for the $k$-th channel (filter) $1 \le k\le c$.
We may slide indices $(i,j)$ by slide size $s$
to reduce the size array in $X^+$.
Otherwise, the dimension $n\times n\times c$ increases as filter size $c$ and results in large storage and increases cost.

First, the depth of the output volume is a hyperparameter: it corresponds to the number of filters we  use, each learning to look for something different in the input. For example, if the first convolution layer takes as input the raw image, then different neurons along the depth dimension may activate in presence of various oriented edges, or blobs of color. We will refer to a set of neurons that are all looking at the same region of the input as a depth column (some people also prefer the term fibre).
Second, we must specify the stride with which we slide the filter. When the stride is 1 then we move the filters one pixel at a time. When the stride is 2 (or uncommonly 3 or more, though this is rare in practice) then the filters jump 2 pixels at a time as we slide them around. This will produce smaller output volumes spatially.

\subsection{Maxpooling}
Max pooling is a down-sampling operation and  commonly used feature extraction operation, typically applied in convolutional neural networks. It reduces the spatial dimensions of features by selecting the maximum value within each small window or region. it is a commonly used feature extraction operation, typically applied in convolutional neural networks. It reduces the spatial dimensions of features by selecting the maximum value within each small window.
Max pooling is a type of pooling operation in convolutional neural networks (CNNs) that reduces the spatial dimensions of a feature map by taking the maximum value within each region. It's often used after convolutional layers to downsample the input, making the network computationally efficient and reducing the number of parameters. 
    
Max Pooling is a pooling operation that calculates the maximum value for patches of a feature map, and uses it to create a downsampling.
Max pooling helps reduce the spatial dimensions of the feature map, making the network less sensitive to small variations in the input data. 
A window (usually $2\times2$) slides across the feature map, and for each window, the maximum value is selected. 
Strides:The stride determines how much the window moves at each step. A stride of 2 means the window moves 2 pixels at a time, reducing the output size by half. 

Benefits are: 
\begin{itemize}
    \item Dimensionality reduction: reduces the size of the feature map, leading to fewer parameters and computational costs.
    \item Translation invariance: makes the network less sensitive to small shifts in the input.
    \item Overfitting reduction: Helps prevent the network from memorizing the training data by reducing the number of features.
\end{itemize}

\subsection{Batch normalization}
We discuss \textbf{Batch Normalization} \cite{IS}, a procedure that accelerates and stabilizes training by reducing internal covariance shift. Let $d$ be the amount of pixel. Let $B$ denote a mini-batch of size $m$, let $\{x_1, \cdots, x_m\}$ be a mini-batch of activations on the $j$-th layer ($1\le j \le d$) , so that the mean of batch $B$ is $\mu_B = \frac{1}{m} \sum_{i=1}^m x_i$, and $\sigma_B^2 = \frac{1}{m} \sum_i (x_i -\mu_B)^2.$ For $k \in \{1, \cdots, d\}$, and a given $\epsilon > 0$, we have: 
\begin{equation}\label{BatNormal}
    \hat{x_i} = \frac{x_i - \mu_B}{\sqrt{(\sigma_B)^2 + \epsilon}}. 
\end{equation}
 
Then, one sets internal variables of a Neural Network as offset parameter $\beta$ and scale parameter $\gamma$ so that $y_i = \gamma \hat{x}_i + \beta.$ 
In terms of the target, instead of a vector with target classes, we consider to use a vectorized target: binary classification with SVM, or multi-class classification with softmax function as a distribution. A multi-class classification generalizes binary classification when it has more than 2 classes of data for making classification. 

\section{Stochastic gradient descent (SGD) for NN design}\label{SGD}
\subsection{Stochastic gradient method}
Consider the minimization of the form
$$
\min \sum_{k=1}^{K} f_k(x)
$$
over randomly sampled indices $k$ where $f_k$ are the loss functions.
For example given input-output sample pairs $(A_k,Y_k)$:
$$
f_k=|\psi((W,b),A_k)-Y_k|^2\mbox{  (a sample $k$)}
$$
where $x=(W,b)$ be the design parameters.
In general a sample size $N$ is large and  so as the dimension $m$ of design parameters $x=(W,b)$.
We randomly sort the samples.
The standard gradient method updates solution $\{x_i\}$ is given by
$$
x_{i+1}=x_i-\alpha\ \sum_{k\in I_i}\nabla f_k(x_i)
$$
where $\alpha>0$ is the activation rate and $\nabla f_k$ is the gradient.
We evaluate the gradient by the back projection (adjoint) method.
i.e, we use the forward and the  backward integration of \eqref{DNN}.

\subsection{Mini-batch method and Anderson-type method}
Consider the optimization of the form
$$
\min\quad \sum_{k\in I} f_k(x)
$$
where $I$ is the sample set with size $N$.
Let $I_i$ be a randomly selected (non-overlapped) sub-indices of 
$I=\cup I_i$. We solve 
$$
\min \quad \sum_{k\in I_i} f_k(x) \mbox{ over the block  $I_i$},
$$
by the subspace method
$$
x_{i+1}=x_i-\sum_{k\in I_i} \beta_k f^\prime_k(x_i),
$$
where $\{\beta_k\}$  minimizes
$$
\min\quad \sum_{k\in I_i} f_k(x_i-\sum_{k\in I_i}\beta_k f^\prime_k(x_i)),
$$
in which $x_i$ is hold on $I_i$. Or,
$$
\min\quad|\sum_{k\in I_i} \beta_k f^\prime_k(x_i)|^2 \mbox{  subject to  }\sum\beta_k=1.
$$
If we select  $\beta_k=\beta$, then
$$
\min  \quad \sum_{k\in I_i} f_k(x_i-\beta\,\sum_{k\in I_i}f^\prime_k(x_i))
$$
where
$$
\frac{1}{m}\sum_{k\in I_i} f^\prime_k(x_i)
$$
is the average of gradients over mini-batch $I_i$.
\section{Backpropagation}
We finally launch the \textbf{Backpropagation Algorithm} \cite{RHW}. It extends the \eqref{DNN}, also known as forward pass, with an additional step of internal variables update, commonly referred as training. Given:
\begin{itemize}
    \item Training data $\{(x_{i}, y_{i})\}_{i=1}^{m}$, with $m$ as batch size;
    \item Fixed Learning rate $\eta$;
    \item Initialized Weights $W_{l}$ and biases $b_{l}$ for layer $l$ with $l \in \Z_+\cap \Z_{\le L}$ and $L$ as total number of layers;
    \item Loss Function $\mathcal{L}$ used to train a Neural Network.
\end{itemize}

\noindent For each training example $(x, y)$, perform the following steps:

1. Forward Pass: We implement \eqref{DNN} step-by-step.  
\begin{align*}
x_{1} &= x \\
\text{For } l = 2 \text{ to } L: \quad
z_{l} &= W_{l} x_{l-1} + b_{l} \\
x_{l} &= \phi_l (z_{l}).
\end{align*}

2. Backward Pass: Using chain rule from calculus, we compute the gradient of the loss function $\mathcal{L}$ over mini-batches with respect to each weight $W$ and bias $b$ to perform gradient descent: 
\begin{align*}
\delta_{L} &= \nabla_{x_L} \mathcal{L}(x_{L}, y) \odot \phi'(z_{L}) \\
\text{For } l = L-1 \text{ to } 2: \quad
\delta_{l} &= (W_{l+1})^\top \delta_{l+1} \odot \phi'(z_{l}),
\end{align*}
where $\odot$ denotes the Hadmard product of two matrices with the same dimension. That is, given two matrices \( A = [a_{ij}] \in \mathbb{R}^{m \times n} \) and \( B = [b_{ij}] \in \mathbb{R}^{m \times n} \), their \textbf{Hadamard product} (denoted by \( A \odot B \)) is defined as:
\[
A \odot B = [a_{ij} \cdot b_{ij}] \in \mathbb{R}^{m \times n}.
\]

3. Parameters' Update:
\begin{align*}
\text{For } l = 2 \text{ to } L: \quad
W_{l} &\gets W_{l} - \eta \ \delta_{l} (x_{l-1})^\top \\ 
b_{l} &\gets b_{l} - \eta \ \delta_{l}. 
\end{align*}
\section{Convergence Stochastic gradient decent and Kaczmarz method}

Recall the stochastic gradient decent method
$$
x_{i+1}=x_i-\alpha_i \sum_{k \in {I_i}}f^\prime_k(x_k),
$$
where $I_i$ is a overlapped mini-batch block .
Consider the case when $f_k=|g_k(x)|^2$ and $I_i={i}$
$$
\min \quad \sum_k |g_k(x)|^2,
$$
and the Kaczmarz update by
$$
x^+-x=-\frac{g_k^\prime(x)}{|g_k^\prime(x)|^2}g_k(x)
$$
Thus, Kaczmarz method \cite{SV} can be seen as a particular (stepsize) case of stochastic gradient descent method. Since
$$
g_k(x^k)=g_k(x^k)(1-g^\prime(\eta)-g^\prime)(x_k))\sim 0
$$
where $\eta$ is the intermediate variable for the gradient.
$x_{k+1}$ is an approximated projection of $x_k$ on the manifold $g_k(x)=0$.
It implies that $\{x_k\}$ converges with a linear rate.

In general one can use the parallel form:
we solve the least square problem  over mini-batch block  $I_i$.
$$
|A_i(x-x_{i}))+A_ix_{i}+gx_{i}|^2 \mbox{  over } x
$$
The solution is given by
$$
x_{i+1}=x_{i}+A_i^t(A_iA_i^t)^{-1}g_k(x_i))
$$
We note that $P=A_i^t(A_iA_i^t)^{-1}A_i$ defines the projection on $R^n$, I.e.,
$$
PP=A_i^t(A_iA_i^t)^{-1}A_iA_i^t(A_iA_i^t)^{-1}A_i=P
$$
where $A_i=\vec{g}^\prime(x_i)$  is the block Jacobean. Then, 
$$
g_k(x_i)-g_k(x^*)\sim A_i(x_i-x^*).
$$
Thus,
$$
|x_{i+1}-x^*|\le \rho_k|x_{i}-x^*|
$$
for some random rate $\rho_i<1$.



\section{Description of algorithm}\label{alg}
We summarize previous sections for a deep neural network and a convolution neural network.
\subsection{DNN Algorithm}
In this section description of the proposed algorithm  for DNN and CNN, specifically applied the EMNIST digit classification.
Define epoch randomly sorted samples $A=X_0$
and $\bar{Y}$ by

We first describe the algorithm on DNN. The algorithm consists of the forward path and 
the backward projection as
\begin{itemize}
\item Initialization of internal variables and input data scaling with division over maximum number in the training matrix.
\item Forward path by DNN: We transport $X_0$ to $Y$ by
$$
X_{k+1}=\phi(W_kX_k+b_k),\;\; X_0=A,\;\;Y=H(X_L).
$$
with Relu $\phi: x \to \max(0,x)$ and $H$ is the softmax function defined in \eqref{SoMax}.
\item Backward projection for the gradient with respect to $W$ and $b$ by, with $B$ as batch size:
\begin{align}
\Delta_L &= Y - \bar{Y}, \\
\Delta_k &= (W_{k+1}^\top \Delta_{k+1}) \odot \mathbf{1}_{Z_k > 0};\\
gW_k &= \tfrac{1}{B}\, \Delta_k A_{k-1}^\top, \\[4pt]
gb_k &= \tfrac{1}{B}\, \Delta_k \mathbf{1}.
\end{align}
\item Solution update by the gradient decent method
$$
W_k^+=W_k- \frac{\alpha}{N}\ gW_k,\;\;b_k^+=b_k-\frac{\alpha}{N}\ gb_k
$$
by a sequence of mini-batches.
\end{itemize}

\begin{remark} In our case, we take average with sample size $N$ so that the gradient descent step has become:
$$
W_k^+=W_k-\frac{\alpha}{N}\ gW_k,\;\;b_k^+=b_k-\frac{\alpha}{N}\ gb_k.
$$
\end{remark}
This diagram is universally defined if the map $H$ and data
$(X_0,\bar{Y})$ are given. 

We pick other epoch and restart the gradient algorithm
for the next update for $(W_k,b_k)$. After 4 epochs
we have a sequence of design parameters $(W_k,b_k)$.
Then we use the Anderson-type method \eqref{And} to accelerate to obtain the improved update.
\subsection{CNN Algorithm}
CNN transfers image $A=X_0$ 
to $B$ with a reduced dimension by applying convolution and max-pooling operations. Then we 
connect $B\to Y$ by a single layer DNN.
The forward path and back projection is governed by convolution
and  adjoint convolution. 

\begin{enumerate}
    \item Convolution Layer:

Given input feature map $X \in \mathbf{R}^{H \times W}$ and kernel $K \in \mathbf{R}^{k_H \times k_W}$, the convolution output at location $(i,j)$ is:

\[
Z_{i,j} = (X * K)_{i,j} = \sum_{m=0}^{k_H-1} \sum_{n=0}^{k_W-1} X_{i+m, j+n} \ K_{m,n}.
\]

For multiple channels and filters:

\[
Z^{(l)}_k = \phi^{(l)} \left( \sum_{c=1}^{C} X^{(l-1)}_c * K^{(l)}_{k,c} + b^{(l)}_k \right).
\]
where, for each $1\le l \le L$:
\begin{itemize}
    \item $X^{(l-1)}_c$: input feature map from channel $c$ at layer $l-1$;
    \item $K^{(l)}_{k,c}$: kernel connecting input channel $c$ to output channel $k$ at layer $l$;
    \item $b^{(l)}_k$: bias for the $k$-th filter at layer $l$;
    \item $\phi^{(l)}$: activation function [e.g., ReLU (Rectified Linear Unit)] at layer $l$.
\end{itemize}

\item Activation Function: Apply a commonly used activation function, say RELU:
\begin{equation}\label{RELU}
\phi(x) = \max(0, x) \quad \text{(ReLU)},    
\end{equation}

\item Fully Connected Layer: After flattening the feature maps, for $l \in \Z_+ \cap \Z_{\le L}$, we have:
\[
a^{(l)} = \phi^{(l)}(W^{(l)} a^{(l-1)} + b^{(l)}).
\]
Typically one uses softmax \cite{BrP}\cite{BrT} as the last activation function for classification tasks:  Let $j \in  \Z_{\ge 1} \cap \Z_{\le n_{L}}$, and $i\in \Z_+\cap \Z_{\le N}$, for the $i$-th sample, let $a^{j}_i$ be the $j$-th coordinate, 
\[\hat{y}_i = \frac{1}{\sum_{j} e^{a_i^{j}}}(e^{a_i^{1}}, \cdots, e^{a_i^{n_L}}). \]
\item  Loss Function: For classification with softmax from above:  
\[
\mathcal{L}(y, \hat{y}) = \frac{1}{2N} \sum_i |y_i - \hat{y}_i|^2.
\]
\end{enumerate}
Thus, the optimization is carried over the internal variables 
$(W,b)$.
\subsection{Accuracy Measurement}
We perform testing by doing \textbf{accuracy measurement} by looking at the image of a vector and see the difference between the predicted value and the actual value. If they are the same, we count; we do not count otherwise. 

We report the number of accuracy after finishing one epoch before moving onto the next epoch, and stop as long as accuracy gets 100\%, or meet stop criterion. For example, we may stop when accuracy level gets above $99.8 \%$, or epoch number gets $1000$.

\section{Numerical Results}
In this section we present numerical findings for DNN and CNN for MNIST digit classification problem.
\subsection{Algorithm implementation}
In this section, we discuss how we implement algorithms through Anderson acceleration. So, the main focus is on neural network architecture, and initialization of internal variables, residual storage for Anderson acceleration matrix. 

We discuss dimensions settings as an essential part to the architecture of Neural Network: For layer design, we let $n_0 = 28 \times 28 = 784$ be the image pixel size, let $n_1 = 128$ be the hidden layer size, let $n_2 = 10$ be classification layer size. Now, we discuss the sample and mini batch amount: Let $N = 5000$ be sample size by picking the first $5000$ from MNIST training set; Let $m = 64$ be the batch size. For the variable step update, we fix the learning rate $\alpha = 0.01$.

We now talk about data initialization: We regularize the input data by making it a matrix with all entries in between 0 and 1. Since for all $(i,j)$ we have $\max_{i,j} [X_{ij}] = 255$ and $[X_{ij}] \ge 0$, we let $X_0 = \frac{1}{255} X_0$ be our initialization. Afterwards, we switch attention internal variables. We put weight parameters $(W_1, W_2)$ as random matrices of comfortable dimension, and set $(b_1,b_2)$ as zeros matrices of comfortable dimension. For the discussion of dimension set up, the beginning of the paper gives dimension for $(W_1,W_2)$ and $(b_1, b_2)$. We use mini-batch gradient descent to train our Neural Network. 

We discuss residual storage for self testing. For the $k$-th epoch from testing (same as training), we use $R_k(:) = P_k(:) - Y(:)$ with vectorization as residual to be collected for Anderson acceleration, where $P_k$ is the predicted probability distribution, and $Y$ is the true probability distribution for testing (same as training) epoch $k$. We apply Anderson acceleration into our training and see that it improves prediction accuracy very much.

For self testing, we do \textbf{restarted Anderson}: We do Anderson acceleration after collecting residual from 4 epochs, and then update parameters by minimizing total sum of residuals before the restart. To accelerate, we repeat the entire process for $m$ times. Afterwards, to accelerate further, one repeats the above process 4 times. In tables involving self testing, we record $4 \times m$ as epoch number. We also perform this algorithm on CNN in the next section.








\subsection{DNN} 

We discuss dimensions settings as an essential part to the architecture of Neural Network: For layer design, we let $n_0 = 28 \times 28 = 784$ be the image pixel size, let $n_1 = 128$ be the hidden layer size, let $n_2 = 10$ be classification layer size. Now, we discuss the sample and mini batch amount: Let $N = 5000$ be sample size by picking the first $5000$ from MNIST training set; Let $m = 64$ be the batch size. For the variable step update, we fix the learning rate $\alpha = 0.01$.

We now talk about data initialization: We regularize the input data by making it a matrix with all entries in between 0 and 1. Since for all $(i,j)$ we have $\max_{i,j} [X_{ij}] = 255$ and $[X_{ij}] \ge 0$, we let $X_0 = \frac{1}{255} X_0$ be our initialization. Afterwards, we switch attention internal variables. We put weight parameters $(W_1, W_2)$ as random matrices of comfortable dimension, and set $(b_1,b_2)$ as zeros matrices of comfortable dimension. For the discussion of dimension set up, section \ref{intro} gives dimension for $(W_1,W_2)$ and $(b_1, b_2)$. We use mini-batch gradient descent to train our Neural Network. 

We discuss residual storage for self testing. For the $k$-th epoch from testing (same as training), we use $R_k(:) = P_k(:) - Y(:)$ with vectorization as residual to be collected for Anderson acceleration, where $P_k$ is the predicted probability distribution, and $Y$ is the true probability distribution for testing (same as training) epoch $k$. We apply Anderson acceleration into our training and see that it improves prediction accuracy very much. We use 4 epochs to generate 
DNN designs sequentially  and then apply the Anderson-type acceleration \eqref{And}. In the following table, we show the accuracy of correctly identifying a class for training data as well as for test data.

\begin{table}[h]
    \centering
    \begin{tabular}{|c|c|c|c|}
    \hline
         Epoch & Standard & Anderson \\ 
         \hline
         4 & 57.94\% & 60.14\% \\ 
         \hline
         8 & 60.74\% & 81.04\% \\ 
         \hline
         20 & 86.52\% & 91.12\% \\
         \hline
         40 & 87.00\% & 92.28\% \\
         \hline
         100 & 91.28\% & 94.02\% \\
         \hline
         500 & 97.86\% & 99.52\% \\
         \hline
         1000 & 99.64\% & 99.46\% \\ 
         \hline
    \end{tabular}
    \caption{DNN Self Test: Accuracy Comparison.} 
    \label{DNNSeTe}
\end{table}

We run DNN with 1000 epochs and it achieves 99.64 \% accuracy.
But, the Anderson method requires much less epochs to 
achieve the same accuracy.







\noindent \textbf{Summary}: Table \ref{DNNSeTe} explains the performance of Standard vs Anderson acceleration with self testing data with Anderson improves a lot. Anderson acceleration has robust performance compared to the standard ones. We include more details here:

1. DNN converges slowly to the desirable accurate solution, e.g. in the standard case, achieving accuracy level over $99\%$ requires 1000 epochs. 

2. Anderson acceleration method converges monotonically with faster convergence, e.g. with 500 epochs, we achieve near 100\% accuracy.

3. DNN can be applied to more general cases. For example, we perform Anderson acceleration on CNN with the full training set as self test, as shown in Table \ref{CNNSe}. 

In conclusion, Anderson acceleration works what we predicted, i.e. Anderson can be used to accelerate DNN, as what we assumed, in general. 
\begin{remark} We use Anderson acceleration independent of the learning net, i.e. residual for each design parameters. We use it on a convolution neural network in Section \ref{CNN}.
\end{remark}
\begin{remark}
     When we consider taking mini-batch and learning from it to make predictions, Anderson accelerates inner loop convergence. 
\end{remark} 
\subsection{CNN}\label{CNN}

We describe variables in our convolution neural network:
\begin{itemize}
    \item Convolution2dLayer: $W_1: 3 \times 3 \times 1 \times 8; b_1: 1 \times 1 \times 8$.
    \item Batchnormalization: After getting $\hat{x}$ from \eqref{BatNormal}, one sets the Offset parameter $\beta$ and the Scale parameter $\gamma$ so that $y = \gamma \hat{x} + \beta.$ Here $\gamma: 1 \times 1 \times 8$, and $\beta: 1 \times 1 \times 8$. 
    \item FullyconnectedLayer: $W_2: 10 \times 6272 ; b_2: 10 \times 1$, where $6272 = 8 \times 784.$
\end{itemize}

For self test and testing, we test all training set from the original MNIST. For testing, we load the original MNIST dataset \cite{LC}: a total training sample as 60000, with testing sample 10000. We set 64 as size for each mini-batch with learning rate $\alpha =0.01$. We do stochastic gradient descent to optimize internal variables. 

We discuss residual storage for testing. For the $k$-th epoch from testing, we use $R_k(:) = P_k(:) - Y(:)$ as residual with vectorization to be collected for Anderson acceleration, where $P_k$ is the predicted probability distribution of epoch $k$, and $Y$ is the true probability distribution. We apply Anderson acceleration into our training and see that it improves prediction accuracy very much.

For self testing and testing, we do \textbf{restarted Anderson}: We do Anderson acceleration after collecting residual from 4 epochs, and then update parameters by minimizing total sum of residuals before the restart, we repeat the process $m$ times. In both cases, we record $4 \times m$ as epoch number.

\begin{table}[h]
    \centering
    \begin{tabular}{|c|c|c|}
    \hline
       Epoch & Standard & Anderson \\
    \hline
        $4$ & $99.31\%$ & $99.41\%$ \\
    \hline
        $8$ & $99.83\%$ & $99.83\%$ \\
    \hline
    \end{tabular}
    \caption{CNN Self Test: Accuracy Comparison.} 
    \label{CNNSe}
\end{table}

\begin{table}[h]
    \centering
    \begin{tabular}{|c|c|c|c|}
        \hline
         Epoch & Standard & Anderson \\
         \hline
         4 & 97.28\% & 97.69\% \\ 
         \hline
         8 & 97.87\% & 98.02\% \\ 
         \hline
    \end{tabular}
    \caption{CNN Test: Accuracy Comparison.}
    \label{CNNTe}
\end{table}
Tables \ref{CNNSe} and \ref{CNNTe} suggest that Anderson acceleration has better performance: The improvement is remarkable after 2 times of applying Anderson acceleration. Since we require way less epochs to achieve accuracy, CNN will be better design for Neural Network.

\section{Conclusion}
We bring more observations to conclude the paper. First we compare the performance without the use of acceleration:  
\begin{enumerate}
    \item The use of CNN and DNN to train a Neural network exhibits the use of universal approximation property.
    \item CNN tables suggest that using a rich architecture on a Neural Network may provide a more accurate prediction result, despite CNN requires a high storage demand has increased by the architecture of CNN, e.g. filter number amplifies the data structures significantly, as we can see from the size of weight matrix for the fully connected layer. DNN uses matrix as data, while CNN uses two dimension image as data. Therefore, DNN does quite good which connect input to output, while cnn transfers images to digit classification directly.
    \item The performance is much better for CNN than DNN. Compared to DNN, CNN converges with $100\%$ with much fewer number of epochs.
    \item CNN performs better in terms of accuracy in fewer epochs. DNN per epoch cheaper than CNN. For self testing with all MNIST training samples (60000) as testing set in standard case, operation time per epoch is less for DNN than CNN. DNN: 400 epochs finish in 206.29 seconds; CNN: 8 epochs finish in 178.773 seconds. Therefore, in terms of speed, we have DNN: 0.516 second/epoch; CNN: 22.34 seconds/epoch.

    \item Variants of design can be tried, one can, for example, redefine the Neural Network by including a maxpooling layer. our test indicates that the simplest CNN works quite well with respect to cost and CPU time.
\end{enumerate}
Now, we bring aspects from the use of Anderson acceleration. 
\begin{enumerate}
    \item Anderson does good job in the convergence of DNN and CNN. We have used mini-batch gradient (on a DNN) descent and stochastic gradient descent (on a CNN) with Anderson acceleration, and Anderson makes things much improved. 
    \item The performance is much better for CNN than DNN; So, CNN achieves the level of accuracy with much less CPU time.
    \item Anderson performs very well and is much effective for DNN case. Because of richness of CNN architecture improves Anderson less, but does improve. 
\end{enumerate}

DNN (Deep Neural Network) approximation refers to using deep neural networks to approximate complex functions or models. This involves training a DNN to learn a mapping between inputs and outputs, allowing it to make predictions or decisions similar to a more computationally expensive model. DNNs are often used for tasks like nonlinear model prediction and function approximation.

\section*{Acknowledgments}
T. Xue would like to thank Dr. Bangti Jin and Dr. Fuqun Han in giving him suggestions to revise the paper. This is a refined version from part of T. Xue's Ph.D. thesis, supported by K. Ito.


\end{document}